\documentclass[a4paper,10pt]{amsart}

\usepackage{amsfonts,amsthm,amssymb,amsmath,amscd}

\newcommand{\NN}{\mathbb{N}}
\newcommand{\ZZ}{\mathbb{Z}}
\newcommand{\QQ}{\mathbb{Q}}

\newcommand{\CC}{\mathbb{C}}
\newcommand{\PP}{\mathbb{P}}

\newcommand{\U}{\mathrm{U}}
\newcommand{\E}{\mathrm{E}}

\begin{document}
\title{Elliptic K3 surfaces admitting a Shioda-Inose structure} 
\date{}
\author{Kenji Koike, Yamanashi University}  
\maketitle
\section{Introduction}
An automorphism of a K3 surface $X$ is called symplectic if it acts on $H^{2,0}(X)$ trivially. 
Such automorphisms were studied by Nikulin in \cite{N1}. He proved that a symplectic involution $\iota$ has eight fixed points and 
the minimal resolution $Y \rightarrow X / \left< \iota \right>$ of eight nodes is again a K3 surface. 
In \cite{SI}, Shioda and Inose proved that every K3 surface $X$ with maximal Picard number $20$ 
has a symplectic involution $\iota$  such that $Y$ is a Kummer surface, and that 
the rational quotient map $\pi : X \dashrightarrow Y$ induces a Hodge isometry $T_X(2) \cong T_Y$, 
where $T_X$ is the transcendental lattice of $X$. 
In general, we say that a K3 surface $X$ admits a Shioda-Inose structure if $X$ has such an involution. 
This definition is due to Morrison (\cite{Mo}), and he proved that 
a K3 surface $X$ admits a Shioda-Inose structure if and only if there exists an Abelian surface $A$ and 
a Hodge isometry $T_X \cong T_A$.
Since the transcendental lattice of a general $(1,d)$-polarized Abelian surface is 
$M_d = \U \oplus \U \oplus \left< -2d \right>$, a K3 surface $X$ with Picard number $17$ admits a Shioda-Inose 
structure if and only if $T_X \cong M_d$, namely $NS(X) \cong \E_8 \oplus \E_8 \oplus \left< 2d \right>$. 
However, to the best of author's knowledge, an explicit example of a $3$-dimensional family of such K3 surfaces with 
the involution $\iota$ is known only for $d=1$ (Appendix in \cite{GL}, \cite{K}) and for $d=2$ (\cite{vGS}). In \cite{K} and \cite{vGS}, 
K3 surfaces $X$ are given as elliptic surfaces with a $2$-torsion section $\sigma$, and $\iota$ is given 
by the fiberwise translation by $\sigma$. 
In this situation, the rational quotient map $X \dashrightarrow Y$ is just an isogeny of degree $2$ between 
elliptic curves over $\CC(t)$, and we have a rational map $Y \dashrightarrow X$ of degree $2$ as the dual isogeny. 
This gives a geometric realization of Kummer sandwich theorem $Y \dashrightarrow X \dashrightarrow Y$ which was proved by Ma (\cite{Ma}).
\\ \indent
In this short note, we show that such pairs of elliptic K3 surfaces exist only for $d = 1,2,3,5,7$ under 
the hypothesis that the Mordell-Weil rank is $0$ (Theorem \ref{Th1}), and we construct $X$ and $Y$ explicitly 
for these values of $d$. 
\section{Elliptic K3 surfaces with a $2$-torsion}
\subsection{}
Let $f : X \rightarrow \PP^1$ be an elliptic K3 surface with the zero section $o$. If $X$ has a $2$-torsion section $\sigma$, 
it is given by the Weierstrass equation
\[
 y^2 = x(x^2 + a(t) x + b(t)), \qquad \deg a(t) \leq 4, \quad \deg b(t) \leq 8
\]
with the projection $f(x,y,t) = t \in \PP^1$, and $\sigma = \{ x=y=0 \}$. Let $\iota$ be the translation by $\sigma$.
It is a Nikulin involution, and we have a K3 surface $Y$ by resolving eight nodes on $X / \left< \iota \right>$. 
The rational quotient map $\phi : X \dashrightarrow Y$ is regarded as an isogeny between elliptic curves over $K=\CC(t)$ 
with the kernel $\{ o, \ \sigma\}$. The Weierstrass model of $Y$ is
\[
 Y : \quad y^2 = x(x^2 - 2a(t)x + a(t)^2 -4b(t)),
\]
and the isogeny $\phi$ and the dual isogeny $\hat{\phi}$ is given by
\begin{align*}
 &\phi : X \longrightarrow Y, \quad (x,y) \mapsto (\frac{y^2}{x^2}, \frac{y(x^2 - b(t))}{x^2}), \\
 &\hat{\phi} : Y \longrightarrow X, \quad (x,y) \mapsto 
(\frac{y^2}{4x^2}, \frac{y(x^2 - a(t)^2 + 4b(t))}{8x^2})
\end{align*}
(\cite{ST}, Chapter III. 4). We denote the projection $(x,y,t) \mapsto t$ by $g : Y \rightarrow \PP^1$.
Up to constants, the discriminants of $X$ and $Y$ are
\[
 \Delta_X(t) = b^2 (a^2 - 4 b), \qquad \Delta_Y(t) = b(a^2 - 4 b)^2.
\] 
For general $a(t)$ and $b(t)$, singular fibers of $X$ and $Y$ are $8 I_1 + 8 I_2$
 and Mordell-Weil groups are $X(K) \cong Y(K) \cong \ZZ / 2 \ZZ$. Transcendental lattices and 
N\'eron-Severi groups are
\[
 T_X \cong T_Y \cong \U \oplus \U \oplus N, \quad
NS(X) \cong NS(Y) \cong \U \oplus N
\]
where $N$ is the Nikulin lattice (\cite{vGS}). 
\subsection{}
We are interested in $a(t)$ and $b(t)$ such that the transcendental lattice $T_X$ of the corresponding K3 surface $X$ 
is $M_d = \U \oplus \U \oplus \left< -2d \right>$. To find such $a(t)$ and $b(t)$, let us study configurations of 
possible singular fibers. For our purpose, Shimada's list (\cite{Shim} and \cite{BK}) is useful, but here 
we make arguments self-contained as possible.  
We denote the simple points of a singular fiber $f^{-1}(\nu)$ by $f^{-1}(\nu)^{\sharp}$, 
which has a natural group structure.
Since the specialization map $X(K)_{tor} \rightarrow f^{-1}(\nu)^{\sharp}_{tor}$ on the torsion subgroup 
is injective (\cite{Mi}, Corollary VII.3.3) and $\sigma \in X(K)$ is of order two, $X$ admits singular fibers of Kodaira's type 
$I_n$, $I_n^*$, $III$ and $III^*$. 
Fundamental invariants for these fibers are summarized in the following table, where 
$L_{\nu}$ is the (negative definite) Dynkin lattice generated by components which do not intersect with $o$, 
$m_{\nu}$ is the number of components, $m^{(1)}_{\nu}$ is the number of simple components, 
$n_{\nu}$ is the number of fixed points by $\iota$ and $c(t) = a(t)^2 - 4 b(t)$.
\begin{table}[htb] \begin{center} \begin{tabular}{|c||c|c|c|c|c|c|c|}
\hline
$f^{-1}(\nu)$ & \multicolumn{2}{|c|}{$I_n$} & \multicolumn{2}{|c|}{$I_{2k}^*$} & $I_{2k+1}^*$ & $III$ & $III^*$ 
\\ \hline
$f^{-1}(\nu)^{\sharp}$ & \multicolumn{2}{|c|}{$\CC^* \times (\ZZ / n \ZZ)$} & \multicolumn{2}{|c|}{$\CC \times (\ZZ / 2 \ZZ)^2$} & 
$\CC \times (\ZZ / 4 \ZZ)$ & \multicolumn{2}{|c|}{$\CC \times (\ZZ / 2 \ZZ)$}
\\ \hline
$L_{\nu}$ & \multicolumn{2}{|c|}{$\mathrm{A}_{n-1}$} & \multicolumn{2}{|c|}{$\mathrm{D}_{2k+4}$} & 
$\mathrm{D}_{2k+5}$ & $\mathrm{A}_1$ & $\E_7$ 
\\ \hline
$\mathrm{ord}_{\nu} \Delta_X(t)$ & \multicolumn{2}{|c|}{$n$} & \multicolumn{2}{|c|}{$2k+6$} & $2k + 7$ & $3$ & $9$
\\ \hline
$m_{\nu}(X)$ & \multicolumn{2}{|c|}{$n$} & \multicolumn{2}{|c|}{$2k+5$} & $2k+6$ & $2$ & $8$
\\ \hline
$m_{\nu}^{(1)}(X)$ & \multicolumn{2}{|c|}{$n$} & \multicolumn{2}{|c|}{$4$} & $4$ & $2$ & $2$
\\ \hline
$\iota$ & (i) & (ii) & (i) & (ii) & (i) & - & -
\\ \hline 
$n_{\nu}(X)$ & $n$ & $0$ & $2k+2$ & $2$ & $2k+3$ & $1$ & $3$
\\ \hline
$g^{-1}(\nu)$ & $I_{2n}$ & $I_{n/2}$ & $I_{4k}^*$ & $I_k^*$ & $I_{4k+2}^*$ & $III$ & $III^*$
\\ \hline
$\mathrm{ord}_{\nu} \Delta_Y(t)$ & $2n$ & $n/2$ & $4k+6$ & $k+6$ & $4k + 8$ & $3$ & $9$
\\ \hline
$\mathrm{ord}_{\nu} b(t)$ & $0$ & $n/2$ & $2$ & $k+2$ & $2$ & $1$ & $3$
\\ \hline
$\mathrm{ord}_{\nu} c(t)$ & $n$ & $0$ & $2k+2$ & $2$ & $2k + 3$ & $1$ & $3$
\\ \hline
\end{tabular} \end{center} \end{table}
\\
These are very well known(see e.g. \cite{Mi}, \cite{SS}, \cite{T}), except perhaps $n_{\nu}$ and the type of the fiber 
$g^{-1}(\nu)$ (Last three columns are determined from $\mathrm{ord}_{\nu} \Delta_X(t)$ and the fiber type of $g^{-1}(\nu)$). 
Here we explain the action of $\iota$ on $I_n$ and $I_n^*$. 
First of all, note that an involution on $\PP^1$ has two fixed points, 
and that intersection numbers are preserved by $\iota$, that is, $D_1 \cdot D_2 = \iota^* D_1 \cdot \iota^* D_2$ for divisors $D_i$.
\subsection{$I_n$-fiber}
\label{In-fiber}
Let $\Theta_k \cong \PP^1 \ (k \in \ZZ / n \ZZ)$ be components of a fiber of type $I_n \ (n >1)$, such that 
$\Theta_0$ intersects with the zero section $o$ and that $\Theta_k \cdot \Theta_{k+1} = 1$ (or $2$ if $n=2$.). 
Then we can identify simple points of $\Theta_k$ with $\CC^* \times \{ k \} \subset \CC^* \times \ZZ/n\ZZ$, replacing 
$\Theta_k$ by $\Theta_{-k}$ if necessary. 
There are two possibilities.
\begin{itemize}
\item[(i)] If $\sigma$ intersects with $\Theta_0$ at the point corresponding to $(-1, 0) \in \CC^* \times \ZZ/n\ZZ$, 
then $\iota$ acts on each $\Theta_k$ as an involution and fixed points are $n$ intersection points $\Theta_k \cap \Theta_{k+1}$.
\item[(ii)] If $n = 2m$ and $\sigma$ intersects with $\Theta_m$ at the point corresponding to 
$(\pm 1, m) \in \CC^* \times \ZZ/2m\ZZ$, then $\iota$ switches $\Theta_k$ and $\Theta_{k+m}$ and there is no fixed point. In this 
case, we define a $\QQ$-divisor
\[
 \vartheta_{2m} = \frac{1}{2m} \sum_{k=1}^{2m-1} k \Theta_k \in NS(X) \otimes \QQ.
\]
Note that $\vartheta_{2m} \cdot \Theta_k \in \ZZ, \ \vartheta_{2m} \cdot o = 0, \ \vartheta_{2m} \cdot \sigma = \frac{1}{2}$ and 
\begin{align*}
\vartheta_{2m} \cdot \vartheta_{2m} &= \frac{1}{4m^2} \{ \sum_{k=1}^{2m-1} k^2 \Theta_k \cdot \Theta_k + 
2 \sum_{k=1}^{2m-2} k(k+1) \Theta_k \cdot \Theta_{k+1} \} \\
&= \frac{1}{4m^2}\{-2 \sum_{k=1}^{2m-1} k^2 + 2 \sum_{k=1}^{2m-2} (k^2+k)\} 
= -1 + \frac{1}{2m}.
\end{align*}
\end{itemize}
We shall use $\vartheta_{2m}$ later, to determine the discriminant group $NS(X)^* / NS(X)$.
\subsection{$I_n^*$-fiber}
\label{In*-fiber}
Next, let $\Theta_0, \cdots, \Theta_3$ be simple components, and $\Gamma_1, \cdots, \Gamma_n$ be double components
of a fiber of type $I_n^*$ as in the following figure,
and let $\Theta_0$ be the component which intersects with $o$. 
\begin{figure}[htbp] \begin{center}
\setlength{\unitlength}{1mm}
\begin{picture}(120,22)
\thicklines
\put(10,1){\line(-1,2){10}} \put(3,15){$\Theta_0$} 
\put(20,1){\line(-1,2){10}} \put(13,15){$\Theta_1$}
\put(5,3){\line(5,1){25}} \put(23,8){\line(4,-1){25}}
\put(36,10){double components}
\multiput(47,4)(1,0){15}{\line(1,0){0.5}}
\put(57,0){\line(4,1){22}} \put(70,5){\line(5,-1){25}}
\put(82,0){\line(1,2){10}} \put(92,15){$\Theta_{2}$}
\put(91,-2){\line(1,2){10}} \put(102,15){$\Theta_{3}$}
\put(0,1){$\Gamma_0$} \put(30,1){$\Gamma_1$} \put(65,-1){$\Gamma_{n-1}$} \put(96,-2){$\Gamma_n$} 
\end{picture}
\end{center} 
\end{figure}
\begin{itemize}
\item[(i)] If $\sigma$ intersects with $\Theta_1$, 
then $\iota$ switches $\Theta_0$ and $\Theta_1$, acts on each $\Gamma_k$ and switches $\Theta_2$ and $\Theta_3$. 
In this case, we have $n$ fixed points $\Gamma_k \cap \Gamma_{k+1}$ and another fixed point on $\Gamma_1$ and
on $\Gamma_n$.
\item[(ii)] If $n = 2m$ and $\sigma$ intersects with $\Theta_2$ or $\Theta_3$, then $\iota$ switches 
$\Theta_0 + \Theta_1$ and $\Theta_2 + \Theta_3$, acts on $\Gamma_m$ and switches $\Gamma_k$ and $\Gamma_{2m-k}$. 
In this case, we have $2$ fixed points on $\Gamma_m$. 
\end{itemize}  
\subsection{Lemma} \label{L}
Let $X$ be an elliptic K3 surface with a $2$-torsion and the Mordell-Weil rank $0$.
\\
(1) If the discriminant group $T_X^* / T_X$ has a subgroup $\ZZ / p^e \ZZ$ for an odd prime number $p$,
then $X$ has a $I_n$-fiber with $n = k p^e$ for some $k \in \NN$. 
\\
(2) If the discriminant group $T_X^* / T_X$ has a subgroup $\ZZ / 2^e \ZZ$ with $e \geq 3$, 
then $X$ has a $I_n$-fiber with $n = 2^e k$ for some $k \in \NN$.
\\
(3) If the Picard number $\rho(X)$ is $17$, then a singular fiber of $X$ is one of the following:
\[
 I_1, \cdots, I_8, I_{10}, I_{12}, I_{14}, I_{16}, I_0^*, \cdots, I_6^*, I_8^*, I_{10}^*, III, III^*,
\]
where $I_{10}, I_{12}, I_{14}, I_{16}, I_8^*$ and $I_{10}^*$ are of type (ii). In particular, possible cyclic subgroups of $T_X^* / T_X$
of order $p^e$ are
\[
 \ZZ / 2^e \ZZ \ (1 \leq e \leq 4), \quad \ZZ / 3 \ZZ, \quad \ZZ / 5 \ZZ, \quad \ZZ / 7 \ZZ. 
\]
\\
{\bf Proof.}
Let $L_X \subset NS(X)$ be the sublattic generated by the zero section $o$, a general fiber and 
components of singular fibers which do not intersect with $o$. Then $L_X$ is of finite index in $NS(X)$,  
and we have 
\[
 L_X \subset NS(X) \subset NS(X)^* \subset L_X^*.
\]
Hence $T_X^* / T_X \cong NS(X)^* / NS(X)$ is isomorphic to a quotient of a subgroup of $L_X^* / L_X$. Note that
\[
 L_X \cong \U \bigoplus (\bigoplus_{\Delta(\nu)=0} L_{\nu}), \qquad 
L_X^* / L_X \cong \bigoplus_{\Delta(\nu)=0} L_{\nu}^* / L_{\nu} 
\]
and $L_{\nu}^* / L_{\nu}$ is one of 
\[
 \mathrm{A}_n^* / \mathrm{A}_n \cong \ZZ / n \ZZ, \quad \mathrm{D}_{2k}^* / \mathrm{D}_{2k} \cong (\ZZ / 2 \ZZ)^2, \quad
\mathrm{D}_{2k+1}^* / \mathrm{D}_{2k+1} \cong \ZZ / 4 \ZZ, \quad 
\E_7^* / \E_7 \cong \ZZ / 2 \ZZ
\]
according to $I_n (III)$, $I_{2k}^*$, $I_{2k+1}^*$ and $III^*$. Therefore subgroups $\ZZ /p^e \ZZ$ stated in (1) and (2) 
come from $I_n$-fibers. 
\\ \indent
By the Shioda-Tate formula (\cite{Mi}, Corollary VII.2.4)
\[
 \rho(X) = 2 + \mathrm{rank} \medspace X(K) + \sum_{\Delta(\nu) = 0} (m_{\nu}(X)-1),
\]
and $\sum n_{\nu} = 8$, we see that a possible singular fiber is in the above list. 
\hfill $\Box$
\subsection{Theorem} \label{Th1}
Let $X$ be an elliptic K3 surface with a $2$-torsion section $\sigma$ which  
gives a Shioda-Inose structure. If $T_X \cong M_d$ and $\mathrm{rank} X(K) = 0$, then 
$d$ is one of $1,\ 2,\ 3,\ 5,\ 7$ or $15$. If $d=15$, the singular fibers of $X$ must be 
$6 I_1 + I_2 + I_6 + I_{10}$ and the singular fibers of $Y$ must be $6I_2 + I_4 + I_3 + I_5$.
(As we shall see later, however, this configuration does not realize K3 surfaces with $T_X \cong M_{15}$.)
\\ \\
{\bf Proof.} 
Since $T_X^* / T_X \cong \ZZ / 2d \ZZ$, we see that a prime factor $p$ of $d$ is $2,3,5$ or $7$, and that  
$p^2 \nmid d$ for $p=3,5,7$. We have also $2^3 \nmid d$ since 
\[
 T_Y^* / T_Y \cong M_d(2)^* / M_d(2) \cong (\ZZ /2 \ZZ)^4 \oplus (\ZZ /4d \ZZ).
\]
Now let $q$ be the maximal prime factor of $d$. 
\\
(1) the case of $q = 7$. Let us show that $d = 7$. 
If $2 | d$, then $Y$ has $I_{7m} + I_{8n}$ and only $(m,n)=(1,1)$ agrees with $\sum (m_{\nu}(Y)-1)= 15$. 
However, $I_7 + I_8$ on $Y$ corresponds to $I_{14} + I_{4}$ or $I_{14} + I_{16}$ on $X$, and both cases contradict 
$\sum (m_{\nu}(X)-1)= 15$. 
Therefore $2$ is not a prime factor of $d$. 
If $3|d$, then $X$ has $I_{3m} + I_{7n}$ and only $(m,n) = (1,1),(1,2),(2,1)$ agree with $\sum (m_{\nu}(X) -1)= 15$. 
However, $I_3 + I_7$ has $10$ fixed points by $\iota$, and this contradicts $\sum n_{\nu}(X) = 8$. 
If $X$ has $I_3 + I_{14}$, then singular fibers of $X$ must be $I_3 + I_{14} + 7I_1$ by the conditions 
$\deg \Delta_X(t) = 24$ and $\sum (m_{\nu}(X)-1)= 15$. This contradicts $\sum n_{\nu}(X) = 8$. We see that also $I_6 + I_7$ is impossible 
since it corresponds $I_{12} + I_{14}$ or $I_{3} + I_{14}$ on $Y$. Therefore $3$ is not a prime factor of $d$. 
By a similar argument, we can show that $5 \nmid d$.
\\
(2) the case of $q=5$. If $2 | d$, then $Y$ has $I_{5m} + I_{8n}$ and only $I_5 + I_8$ agrees with 
$\sum (m_{\nu}(Y)-1)= 15$.  
This configuration is given as a degeneration (confluences of singular fibers)
\[
 8 I_1 \ (b(t) = 0) \ + \ 8 I_2 \ (c(t) = 0) \ \dashrightarrow \
 (I_5 + 3 I_1) + (I_8 + 4 I_2).
\]
of the most general configuration $8I_1 + 8I_2$.  
Under the hypothesis $\sum m_{\nu}(Y)= 15$, we may admit only $I_1 + I_2 \dashrightarrow III$ as extra confluences. 
By Corollary 1.7 in \cite{Shio}, we have
\[
 |\det NS(Y)| = \frac{\prod m_{\nu}^{(1)}(Y)}{|Y(K)_{tor}|^2} \leq \frac{5 \cdot 8 \cdot 2^4}{4} 
< |\det M_{10}(2)|.
\]
Therefore we see that $2$ is not a prime factor of $d$, and we have $d=5$ or $d=15$. If $d=15$, then $X$ has one of 
\[
 I_3 + I_5, \quad  I_6 + I_{10}, \quad I_{3} + I_{10}, \quad I_6 + I_5.
\]
As degenerations of $8I_1 + 8I_2$, these are
\[
 (I_3 + I_5) + 8I_2, \quad 8I_1 + (I_6 + I_{10}), \quad (I_3 + I_5) + (I_6 + I_{10})
\]
where $I_{3} + I_{10}$ and $I_6 + I_5$ correspond to the same degeneration. From this, $\sum (m_{\nu}(X)-1)= 15$ 
and the equality
\[
 \det NS(Y) = -\det M_{15}(2) = -2^5 \det M_{15} = 2^5 \det NS(X),
\]
we see that the singular fibers of $X$ and $Y$ must be the stated form.
\begin{align*}
\begin{array}{c|ccccc} 
X & 8 I_2 + 8 I_1 & \dashrightarrow & (I_6 + I_{10}) + 8 I_1 & \dashrightarrow & (I_6 + I_{10}) + (I_2 + 6 I_1) \\
\hline
Y & 8 I_1 + 8 I_2 & \dashrightarrow & (I_3 + I_5) + 8I_2 & \dashrightarrow & (I_3 + I_5) + (I_4 + 6 I_2) 
\end{array}
\end{align*}
(3) the case of $q=3$. If $2|q$, then $Y$ has $I_3 + I_8$ or $I_6 + I_8$, that is, the singular fibers of $Y$ are obtained 
as a degeneration of one of the following two configurations.
\begin{align*}
\begin{array}{c|ccc} 
X & 8 I_2 + 8 I_1 & \dashrightarrow & (I_6 + 5 I_2) + (I_4 + 4 I_1) \ \text{or} \ 8 I_2 + (I_3 + I_4 + I_1) \\
\hline
Y & 8 I_1 + 8 I_2 & \dashrightarrow & (I_3 + 5 I_1) + (I_8 + 4 I_2) \ \text{or} \ 8 I_1 + (I_6 + I_8 +  I_2)
\end{array}
\end{align*}
By the condition $\sum (m_{\nu}(Y)-1)= 15$, we may admit just one of the following confluences
\[
 4I_k \dashrightarrow 2I_{2k}, \quad 3I_k \dashrightarrow I_{3k} \quad (k=1,2), \quad  2(I_1 + I_2) \dashrightarrow I_0^*,
\]
and $I_1 + I_2 \dashrightarrow III$ if possible. In any cases, we have 
\begin{align*}
\frac{\prod m_{\nu}^{(1)}(Y)}{|Y(K)_{tor}|^2} = |\det NS(Y)| &= |\det T_Y| \\
&= 2^5 |\det T_X| = 2^5 |\det NS(X)| = 2^5 \frac{\prod m_{\nu}^{(1)}(X)}{|X(K)_{tor}|^2}
\end{align*}
and $|X(K)_{tor}| = 2^{\varepsilon}|Y(K)_{tor}|$ with $\varepsilon = 1,0,-1$. Therefore we have an inequality
\[
 \prod m_{\nu}^{(1)}(Y) \leq 2^3 \prod m_{\nu}^{(1)}(X)
\]
and we see easily that this contradicts any case of the considering degenerations.
\\
(4) the case of $q=2$. If $4|d$, then $Y$ must have $I_{16}$. In this case, the singular fibers of $Y$ are 
$I_{16} + 8 I_1$ and the singular fibers of $X$ are $I_8 + 8 I_2$. These K3 surfaces are studied in \cite{vGS}, and 
we have $T_X \cong M_2$ and $T_Y \cong M_2(2)$.
\hfill $\Box$
\section{Examples}
\subsection{}
For a cubic polynomial $P(t)$ and $0 \leq n \leq 8$, we define an elliptic K3 surface $X_d = X(d, P)$ by
\[
 y^2 = x(x^2 + P(t)x + t^d).
\]
Then the quotient surface $Y_d = X_d / \left< \iota \right>$ is
\[
 y^2 = x(x^2 -2P(t)x + P(t)^2 - 4t^d).
\]
The singular fibers of $X_d$ and $Y_d$ for a general $P(t)$ are given in the following table
\begin{table}[htb] \begin{center} \begin{tabular}{|c|c|c|c|c|c|c|c|c|}
\hline
 & $X_0$ & $X_d (1 \leq d \leq 6)$ & $X_7$ & $X_8$ & $Y_0$ & $Y_d (1 \leq d \leq 6)$ & $Y_7$ & $Y_8$ 
\\ \hline 
$t=0$ & reg. & $I_{2d}$ & $I_{14}$ & $I_{16}$ & reg. & $I_d$ & $I_7$ & $I_8$ 
\\ \hline
$c(t) = 0$ & $6 I_1$ & $6 I_1$ & $7 I_1$ & $8 I_1$ & $6 I_2$ & $6 I_2$ & $7 I_2$ & $8 I_2$  
\\ \hline
$t = \infty$ & $I_{12}^*$ & $I_{12-2d}^*$ & $III$ & reg. & $I_6^*$ & $I_{6-d}^*$ & $III$ & reg. 
\\ \hline
\end{tabular} \end{center} \end{table}
\\
where $c(t) = P(t)^2 - 4 t^d$.
Elliptic K3 surfaces $X_1$ were studied by Kumar in \cite{K}. The transcendental lattice of a general $X_1$ is $M_1$ and 
the quotient surfaces $Y_1$ are Jacobian Kummer surfaces. Elliptic K3 surfaces $X_8$ were studied by van Geemen and Sarti in \cite{vGS}. 
The transcendental lattice of a general $X_8$ is $M_2$ and the quotient surfaces $Y_8$ have the transcendental lattice $M_2(2)$. 
\subsection{Proposition}
For a general cubic polynomial $P(t)$, we have
\\
(1) the Picard number $\rho(X_d)$ is $17$ for $d=1, \cdots, 6$, and $\rho(X_0) = 18$, \\
(2) $X_d(K) = \{ o, \ \sigma \} \cong \ZZ / 2 \ZZ$ for $d=1, \cdots, 6$, \\
(3) $\det NS(X_d) = 2d$ for $d=1, \cdots, 6$, and $\det NS(X_0) = -1$. Hence $NS(X_0) \cong \E_8 \oplus \E_8 \oplus \U$ 
and $T_{X_0} \cong \U \oplus \U$. \\
(4) $T_{X_d} \cong M_d$ for $d=1, \cdots, 6$.
\\ \\
{\bf Proof.}
(1) Cubic polynomials $P(t)$ form a $4$-dimensional vector space, and we have isomorphisms  
\[
 X(d, \lambda^{-4d}P(\lambda^8 t)) \longrightarrow X(d, P(t)), \quad (x,y,t) \mapsto 
(\lambda^{4d}x, \lambda^{6d}y, \lambda^8 t)
\]
by $\lambda \in \CC^*$. Up to this $\CC^*$-action, the configuration of singular fibers is determined by $P(t)$ and it
gives the moduli of $X(d, P)$ for $d=1, \cdots, 7$. Therefore K3 surfaces $X(d,P)$ form a $3$-dimensional family in this case. 
For $d = 0$, we can transform $P(t)$ into $t^3 + at + b$ by a transformation $t \mapsto \alpha t + \beta$, 
and $(a, b)$ gives the moduli. From this, we see that $\rho(X_d) \leq 17$ for $d=1,\cdots,7$ and $\rho(X_0) \leq 18$. 
On the other hand, by the formula 
\[
 \rho(X_d) = 2 + \mathrm{rank} \medspace X_d(K) + \sum_{\Delta(\nu) = 0} (m_{\nu}(X_d)-1),
\]
we have
\[
 \rho(X_d) = \begin{cases} 18 + \mathrm{rank} \medspace X_0(K) \quad (d=0) \\
17 + \mathrm{rank} \medspace X_d(K) \quad (d=1,\cdots,6) \\ 
16 + \mathrm{rank} \medspace X_7(K) \quad (d=7). \end{cases}
\]
Therefore we see that $\mathrm{rank} \medspace X_d(K) = 0$ for $d=0,\cdots,6$. 
\\
(2) We have an injective homomorphism $X_d(K)_{tor} \rightarrow (\ZZ / 2 \ZZ)^2$ 
since $f^{-1}(\infty)^{\sharp} \cong \CC \times (\ZZ / 2 \ZZ)^2$.
The $2$-torsion subgroup of $X_d$ is given by $o, \ \sigma$ and two solutions of $F(x) = x^2 +P(t) x + t^d = 0$. 
Since $F(x)$ is irreducible over $K$, we have 
\[
 X_d(K) = X_d(K)_{tor} \cong \ZZ / 2 \ZZ
\]
for $d=0,\cdots,6$.
\\
(3) By Corollary 1.7 in \cite{Shio}, we have
\[
 |\det NS(X_d)| = \frac{\prod m_{\nu}^{(1)}(X_d)}{|X_d(K)_{tor}|^2} = \frac{1}{4} \prod m_{\nu}^{(1)}(X_d),
\] 
and $\det NS(X_d) = 2d$ for $d=1, \cdots, 6$, and $\det NS(X_0) = -1$.
\\
(4) The N\'eron-Severi group $N = NS(X_d)$ is generated by $o, \ \sigma$ and all components of singular fibers. 
Since the singular fiber at $t=0$ is an $I_{2d}$-fiber of type (ii), we can define $\vartheta_{2d} \in N \otimes \QQ$ 
as in \ref{In-fiber}.
Let $\Theta_0, \cdots, \Theta_3$ be simple components of $I_{12-6d}^*$-fiber at $t = \infty$ as in \ref{In*-fiber}.
This fiber is of type (ii), and $\sigma$ intersects with either $\Theta_2$ or $\Theta_3$. 
Let us consider 
\[
 \Gamma = \frac{1}{2}(\Theta_2 + \Theta_3) + \vartheta_{2d} \in N \otimes \QQ.
\]
Since the intersection numbers of $\Gamma$ with $o, \ \sigma$ and components of singular fibers are integers, 
we see that $\Gamma \in N^*$. Moreover we have $2d \Gamma \in N$ and the value of the discriminant form 
$q_N : N^* / N \rightarrow \QQ / 2 \ZZ$ for $\Gamma$ is 
\[
\Gamma \cdot \Gamma = \frac{1}{4}\{ (\Theta_2)^2 + (\Theta_3)^2 \} + (\vartheta_{2d})^2  = -2 + \frac{1}{2d} 
\equiv \frac{1}{2d} \mod 2.
\] 
If $m \Gamma \in N$, then we have $(\Gamma, m \Gamma) \in N^* \times N$ and 
\[
 \frac{m}{2d} \equiv \Gamma \cdot (m \Gamma) \equiv 0 \mod \ZZ.
\]
Therefore $\Gamma$ gives an element of order $2d$ in $N^*/N$, and we have $N^*/N \cong \ZZ / 2d \ZZ$. 
By Corollary 1.13.3 in \cite{N2}, we see that $N \cong \E_8 \oplus \E_8 \oplus \left< 2d \right>$ and 
$T_{X_d} \cong M_d$.
\hfill $\Box$
\subsection{Lemma}
\label{Lem}
For a general cubic polynomial $P(t)$, we have
\[
 \det NS(Y_d) = \begin{cases} 2^4 \quad (d=0) \\ 2^6 \cdot d \quad (d=1,3,5) \\ 
2^4 \cdot d \quad (d=2,4,6) \end{cases}.
\]
\\
{\bf Proof.}
Since $Y_d(K)$ is isogeneous to $X_d(K)$, we see that $Y_d(K) = Y_d(K)_{tor}$.
For $d=0,2,4,6$, the group structure at $t=\infty$ is $\CC \times (\ZZ / 2 \ZZ)^2$, and we have 
full two-torsions:
\[
 y^2 = x(x^2 - 2P(t)x + P(t)^2 -4t^d) = x(x -P(t) + 2t^{d/2})(x -P(t) - 2t^{d/2}).
\]  
Therefore 
the Mordell-Weil group $Y_d(K)$ is isomorphic to $(\ZZ/2\ZZ)^2$ in this case. For $d = 1,3,5$, the group structure 
at $t=\infty$ is $\CC \times (\ZZ / 4 \ZZ)$, and we have $Y_d(K) = \ZZ / 2 \ZZ$ or $\ZZ / 4 \ZZ$. 
If $\sigma' \in Y_d(K)$ has order four, then $\hat{\phi}(\sigma') \in X_d(K) \cong \ZZ / 2 \ZZ$ has order two.
However, the non-zero element of $X_d(K)$ is pulled back to a double section
\[
 y = 0, \ x^2 - 2P(t)x + P(t)^2 -4t^d = 0
\]
of $Y_d$ by $\hat{\phi}$. Hence we see that $Y_d(K) \cong \ZZ / 2 \ZZ$. As in the case of $X_d$, 
the Lemma follows from Corollary 1.7 in \cite{Shio}.
\hfill $\Box$
\subsection{Proposition}
Let $P(t)$ be a general cubic polynomial.
\\
(1) The rational map $\phi : X_d \dashrightarrow Y_d$ gives a Shioda-Inose structure for $d=0,1,3,5$. 
In particular, $Y_d$ is a Kummer surface with the transcendental lattice $\U(2) \oplus \U(2)$ for $d=0$, and 
$M_d(2)$ for $n=1,3,5$.
\\
(2) The transcendental lattice of $Y_d$ is $\U(2) \oplus \U(2) \oplus \left< -d \right>$ for $d=2,4,6$.
\\ \\
{\bf Proof.}
We have a natural map $\phi_* : T_{X_d} \rightarrow T_{Y_d}$ between transcendental lattices such that 
$\phi_* T_{X_d} \cong T_{X_d}(2)$ (see \cite{SI} and \cite{Mo}). 
\\
(1) By Lemma \ref{Lem}, we see that $\det T_{X_d}(2) = \det T_{Y_d}$. Therefore we have 
$T_{X_d}(2) \cong T_{Y_d}$.
\\
(2) By Lemma \ref{Lem} and the conditions
\[
 \phi_* T_{X_d} \subset T_{Y_d} \subset (T_{Y_d})^* \subset (\phi_* T_{X_d})^*, \qquad 
(\phi_* T_{X_d})^* / \phi_* T_{X_d} \cong (\ZZ / 2 \ZZ)^4 \times (\ZZ / 4d \ZZ),
\]
we see that $(T_{Y_d})^* / T_{Y_d}$ is isomorphic to one of groups 
\[
 (\ZZ / 2 \ZZ)^2 \times (\ZZ / 4n \ZZ), \quad (\ZZ / 2 \ZZ)^3 \times (\ZZ / 2d \ZZ), \quad
(\ZZ / 2 \ZZ)^4 \times (\ZZ / d \ZZ).
\]
Let us consider a sublattice $L$ of $N = NS(Y_d)$ generated by the zero section, a general fiber and components of singular fibers 
which does not intersect with the zero section. Then we have
\[{}
 L \subset N \subset N^* \subset L^*, \qquad L^* / L = (\ZZ / 2 \ZZ)^8 \times (\ZZ / d \ZZ)
\]
since $Y_d \ (d=2,4,6)$ has singular fibers $I_d$, $6I_2$ and $I_{6-d}^*$. Hence $N^* / N$ does not contain an element of order $2d$, 
nor does $(T_{Y_d})^* / T_{Y_d}$. From this, we see that 
\[
 (T_{Y_d})^* / T_{Y_d} \cong (\ZZ / 2 \ZZ)^4 \times (\ZZ / d \ZZ)
\]
and $T_{Y_d} \cong \U(2) \oplus \U(2) \oplus \left< -d \right>$.
\hfill $\Box$
\\
\subsection{}
Let us consider a family of elliptic K3 surfaces
\[
 X_n' : y^2 = x(x^2 + P(t)x + t^n(t-1)^{8-n}), \quad P(t)= 2 t^4 - (8-n)t^3 + a_1 t^2 + a_2 t + a_3 
\]
for $n=5, 7$. A general $X_n'$ has singular fibers $I_{2n}$, $I_{16-2n}$, $I_2$ and $6 I_1$
at $t=0, 1, \infty$ and $P(t)^2 - 4 t^n(t-1)^{8-n} =0$, respectively. A general $Y_n'= X_7' / \left< \iota \right>$ 
has singular fibers $I_n$, $I_{8-n}$, $I_4$ and $6 I_2$ at $t=0, 1, \infty$ and $P(t)^2 - 4 t^n(t-1)^{8-n} =0$, respectively. 
\subsection{Proposition} 
For a general $P(t)$, we have \\
(1) $T_{X_7'} \cong M_7$ and $T_{X_5'} \cong \U \oplus \begin{bmatrix} 2 & 1 \\ 1 & -2 \end{bmatrix} \oplus \left< -6 \right>$, \\
(2) $T_{Y_n'} \cong T_{X_n'}(2)$ for $n=5,7$.
\\
{\bf Proof.}
(1) Since $\sum (m_{\nu}(X_n')-1) = 15$ and 
$a_1,a_2,a_3$ give the moduli parameters, we have $\rho(X_n') = 17$ and $\mathrm{rank} X_n'(K) = 0$.
Then we have
\[
 \det NS(X_n') = \frac{\prod m_{\nu}^{(1)}(X_n')}{|X_n'(K)_{tor}|^2} = \frac{2n(8-n)\cdot 2^2}{|X_n'(K)_{tor}|^2}.
\]
Since $2n(8-n)$ is square-free for $n=5,7$, we see that $X_n'(K) = \{o, \sigma\}$ and 
\[
 \det NS(X_n') = 2n(8-n) = \begin{cases} 30 \quad (n=5) \\ 14 \quad (n=7) \end{cases}
\]
Note that 
singular fibers at $0$ and $1$ are of type (ii) and we have $\vartheta_{2n},\ \vartheta_{16-2n} \in NS(X_n') \otimes \QQ$. 
Since the singular fiber at $\infty$ is of type (i), the component $\Theta_1$ does not intersects with $o$ and $\sigma$. 
Then $\Gamma = \vartheta_{2n} + \vartheta_{16-2n} + \frac{1}{2} \Theta_1$ belongs to $NS(X_n')^*$ and
\begin{align*}
 \Gamma \cdot \Gamma &= (\vartheta_{2n})^2 + (\vartheta_{16-2n})^2 + (\frac{1}{2} \Theta_1)^2 \\
&= (-1 + \frac{1}{2n}) + (-1 + \frac{1}{16-2n}) + (-\frac{1}{2}) =
\begin{cases} -2 - \frac{7}{30} \quad (n=5)\\ -2 + \frac{1}{14} \quad (n=7) \end{cases}.
\end{align*}
From this, we see that $NS(X_7') \cong \E_8 \oplus \E_8 \oplus \left< 14\right>$ and $T_{X_7'} \cong M_7$. 
Let $e_1, \cdots, e_5$ be the basis of 
$M = \U \oplus \begin{bmatrix} 2 & 1 \\ 1 & -2 \end{bmatrix} \oplus \left< -6 \right>$. Then 
we have 
\[
 \delta = \frac{1}{5}(2 e_3 + e_4) + \frac{1}{6} e_5 \in M^*
\]
and $\delta$ generates $M^* / M \cong \ZZ / 30 \ZZ$. since $\delta \cdot \delta = \frac{7}{30}$, we see that 
$T_{X_5'} \cong M$ by Corollary 1.13.3 in \cite{N2}.
\\
(2) We see easily that $Y_n'(K) \cong \ZZ / 2 \ZZ$ and $\det NS(Y_n') = 2^5 \det NS(X_n')$. Hence 
we have $T_{Y_n'} \cong \phi_* T_{X_n'} \cong T_{X_n'}(2)$.
\hfill $\Box$

\end{document}